\documentclass{amsproc}
\usepackage{eucal}
\newcommand {\C}{{\mathbb C}}
\newcommand {\cl}{{\rm cl \;}}
\newcommand {\Diff}{{{\rm Diff}^{+}}}
\newcommand {\G}{{\mathbb G}}
\newcommand {\e}{{\bf e}}
\newcommand {\End}{{\rm End}}
\newcommand {\HH}{{\mathbb H}}
\newcommand {\Hom}{{\rm Hom}}
\newcommand {\kw}{{\rm kw}}                                                                                                                                                                                                                                                                                                                                                                                                                                                                                                                                                                                                                                                     
\newcommand {\mL}{{\mathbb L}}

\newcommand {\Q}{{\mathbb Q}}
\newcommand {\R}{{\mathbb R}}
\newcommand {\SP}{{\rm SP}}
\newcommand {\T}{{\mathbb T}}
\newcommand {\U}{{\mathbb U}}
\newcommand {\V}{{\rm VOA}}
\newcommand {\W}{{\mathbb W}}
\newcommand {\Y}{{\mathbb Y}}
\newcommand {\Z}{{\mathbb Z}}
\newcommand {\MM}{{\overline {\mathcal M}}}
\newcommand {\M}{{\mathcal M}}
\newcommand {\MU}{{\bf MU}}
\newcommand {\X}{{\mathcal Q}}
\newcommand {\half}{{\frac{1}{2}}}
\newcommand {\la}{{\langle}}
\newcommand {\ra}{{\rangle}}

\begin{document}
 
\title{Schur $Q$-functions and a Kontsevich-Witten genus}
 
\author{Jack Morava}
\address{Department of Mathematics, Johns Hopkins University,
Baltimore, Maryland 21218}
\email{jack@math.jhu.edu}
\thanks{The author was supported in part by the NSF}
\date{3 January 1998}
\begin{abstract} The Virasoro operations in Witten's theory of
two-dimensional topological gravity have a homotopy-theoretic
interpretation as endomorphisms of an {\bf ordinary} cohomology 
theory with coefficients in a localization of I.\ Schur's ring 
$\Delta$ of $Q$-functions. The resulting theory has many of the
features of a vertex operator algebra.\end{abstract}

\maketitle

\section*{\bf Introduction}

Smooth complex curves of genus $g>1$ form a stack $\M_{g}$ with
a compactification defined by adjoining the divisor of stable singular 
curves (with double points, but only finitely many automorphisms). Ideas 
from string theory have led to great progress in understanding the topology 
of these moduli objects in the large genus limit. This note is concerned 
with an algebra introduced by Schur in connection with the classification of 
projective representations of symmetric groups, and its relevance to these
spaces. This algebra appears in disguise in classical work on the Riemann 
moduli space and in more recent work of Witten and Kontsevich on the 
intersection theory of its Deligne-Mumford compactification $\MM_{g}$, but 
its significance has become clear only in retrospect, as integral 
cohomology emerges from the fog of $\Q$.

The central construction of topological gravity is a partition function 
which can be defined geometrically by a family $$\tau_{g} : \MM_{g} 
\rightarrow \MU \hat \otimes \Q[v]$$ of maps to the complex cobordism spectrum
tensored with the rationals, and the main result of this paper is the 
construction of a morphism $$\kw : \MU \rightarrow {\bf H}(\Delta[q_{1}^{-1}])
$$ of ring-spectra, which (when composed with $\tau$) sends the fundamental 
class of the moduli space of stable curves to a highest-weight vector for 
a naturally defined Virasoro action on $\Delta_{\Q}$. This Kontsevich-Witten 
genus can be constructed in purely algebraic terms from the theory of formal 
groups, but its natural context is the topology of moduli spaces. The first 
section below is a summary of some background from algebraic geometry, 
emphasizing homotopy theory, while the Hopf algebra of $Q$-functions is the 
topic of the second. This whole subject is in many ways still quite mysterious, 
and a final section argues that some features are most natural in an 
equivariant context. An appendix summarizes the construction of a vertex 
operator algebra following ideas of A. Baker.

I would like to thank M.\ Ando, A.\ Baker, F.\ Cohen, C.\ Y.\ Dong, C.\ Faber, 
T.\ J\'ozefiak, M.\ Karoubi, E.\ Looijenga, G.\ Mason, S.\ Morita, and U.\ Tillmann 
for conversations and correspondence about the material in this paper.

\section*{\bf \S 1 Background from geometry}

\subsection*{1.1} Teichm\"uller theory describes the Riemann moduli space 
$\M_{g}$ as the quotient of a complex $3(g-1)$-dimensional cell by an
action of the mapping class group $$\Gamma_{g} = \pi_{0} \Diff
(\Sigma_{g})$$ of isotopy classes of orientation-preserving diffeomorphisms 
of a closed surface $\Sigma_{g}$; this action is properly discontinuous 
and almost free, in the sense that its isotropy groups are finite, so the 
resulting quotient is an orbifold. The components of the group of 
diffeomorphisms are contractible, and the map $$B \Diff (\Sigma_{g}) = 
B\Gamma_{g} \rightarrow \M_{g}$$ from the homotopy-theoretic quotient of the 
Teichm\"uller action to the geometric quotient is a rational homology 
isomorphism. 

There is no such essentially topological description of
the compactification $\MM_{g}$, but there are very interesting
(proper) forgetful maps $$\Phi^{n}_{g} : \MM^{n}_{g} \rightarrow \MM_{g}$$ 
defined on the moduli stack of stable curves marked with $n$
ordered smooth points. [A marked curve is stable, as above, if it
has only finitely many automorphisms; the loss of the marking
may render an irreducible component of the curve unstable, and
the forgetful map is understood to contract such a component to a
point.] According to the conventions of Quillen, such a proper
complex-oriented map between homology manifolds defines an
element $$[\Phi^{n}_{g}] \in MU^{-2n}_{\Q}(\MM_{g})$$
of the ring (tensored with $\Q$) of complex cobordism classes of maps to the moduli
space. These cobordism classes can be
identified in terms of ordinary cohomology by the ring isomorphism
$$MU^{*}_{\Q}(X) \rightarrow H^{*}(X,\Q[t_{k} | k \geq 1])$$
which sends $[\Phi : V \rightarrow X]$ to the characteristic number polynomial
$$\sum_{I} \Phi_{*}m_{I}(T_{\Phi}^{*})t^{I} \;,$$ where $T_{\Phi}^{*}$ 
is the formal cotangent bundle along $\Phi$ and $\Phi_{*}$
is the Gysin or transfer map in cohomology; $I = i_{1},\dots$ is a
multi-index, $t^{I} = \prod t_{k}^{i_{k}}$, and $m_{I}$ denotes the
characteristic class associated to the monomial symmetric
function indexed by the partition of $|I| = \sum k i_{k}$ 
defined by  $I$. [The $t_{k}$ correspond to the complete symmetric 
functions $h_{k}$ but we keep the standard notation.] This 
differs from the usual [29 \S 6.2] definition, which is
expressed in terms of the normal bundle rather than the cotangent
bundle, but the two constructions are related by a
straightforward change of basis for the characteristic classes,
and the choice above leads most directly to the formulas of
Witten.
   
\subsection*{1.2} The action of a diffeomorphism of the surface $\Sigma_{g}$ on
its homology defines an integral symplectic representation of $\Gamma_{g}$, 
and thus a map $$B \Diff (\Sigma_{g}) \rightarrow BSp(2g,\Z) \; .$$ The group $Sp(2g,\Z)$
acts with finite isotropy on the contractible
symmetric space $Sp(2g,\R)/\U(g)$, with quotient the space $\mathcal A_{g}$ 
of $g$-dimensional principally polarized abelian varieties, so the homomorphism induced by this map 
on rational homology agrees with that defined by the construction which
assigns to a Riemann surface its Jacobian. By stability theorems of Harer and 
Ivanov, these maps are suitably compatible for increasing $g$, and we can think
of them as taking values in the classifying space $BSp(\Z)$ of the
infinite integral symplectic group. The composition of the
obvious maps $$BSp(\Z) \rightarrow BSp(\R) = B\U \rightarrow B(\U/{\mathbb O}) 
= Sp(\HH)/\U$$ with a final Bott isomorphism is a rational
homology isomorphism which (away from two [20 \S 3.15]) splits 
a copy of $Sp(\HH)/\U$ off the group completion $BSp(\Z)^{+}$, cf.\ [5,30]. 

\subsection*{1.3.1} One of the main contentions of this paper is that the
theory of topological gravity defines maps to the complex
cobordism spectrum which are natural analogues, for the compactifications
$\MM_{g}$, of the classical Abel-Jacobi map described above. The physicists' 
definitions are motivated by ideas from statistical mechanics: using the 
language of cobordism, let $$\tau_{g} = \sum_{g \geq 0} [\Phi^{n}_{g}] 
\frac{v^{n}}{n!} \in MU^{0}_{\Q[v]}(\MM_{g}) \; ,$$ with $v$ a bookkeeping variable 
of (cohomological) degree two. The element $\tau_{g}$ is the class of the space 
of configurations defined by an indefinite number of distinct but
unordered smooth points on a stable curve of genus $g$; it is
tempting to think of this as the space of states of the `Mumford
gas' of free particles on a Riemann surface, with existence
as their only attribute. The homotopy class 
representing this tau-function induces a homomorphism $$\tau_{g*} : 
H_{*}(\MM_{g},\Q) \rightarrow H_{*}(\MU,\Q[v]) = \Q[t_{k}|k \geq 1][[v]] 
\; ;$$ the image of the fundamental class $[\MM_{g}]$ of the moduli space 
under this homomorphism
is essentially Witten's free energy $F_{g}$. [It is convenient to 
interpret homology with $\Q[v]$ coefficients to be $v$-adically completed;
for a finite complex this implies no change at all.]

\subsection*{1.3.2} To be more precise about this, 
we will need the slightly more
sensitive characteristic number homomorphism which assigns to
$[\Phi] \in MU^{-2n}(X)$ the class $$\cl [\Phi] = \sum
\Phi_{*}m_{I}(T_{\Phi}^{*}) t_{0}^{n-l(I)}t^{I} \;,$$ where $l(I)
= \sum_{k > 0} i_{k}$ is the length (or number of parts) of the
partition $I$. The ring $\Q[t_{k}|k \geq 0][t_{0}^{-1}]$ of
cohomology coefficients has now been enlarged to include an
invertible polynomial generator of degree zero; this corresponds
[cf. \S 1.5 below] to a natural extension of the Landweber-Novikov
algebra of cobordism operations. The monomial symmetric function 
$m_{I}$ is a sum
$$\sum x_{\sigma(1)}^{d_{1}} \dots x_{\sigma(n)}^{d_{n}}  \;,$$
where $d_{i}$ is a finite sequence of nonnegative numbers with
$\sum d_{i} = |I|$ in which $k$ appears $i_{k}$ times; it is
convenient to think of this sequence as having exactly $n$ terms,
with zero appearing $i_{0} = n - l(I)$ times. The stable
cotangent bundle $T_{\Phi}^{*}$ along $\Phi_{g}^{n}$ is the sum
of the cotangent line bundles $L_{i}$ of the modular curve at its
marked points, and the characteristic class $m_{I}(T_{\Phi}^{*})$
is obtained by substituting Euler classes $\e(L_{i})$ for
the formal variables $x_{i}$. Witten's characteristic number [31
\S 2.4] $$\langle \tau^{I} \rangle = (\prod_{1 \leq i \leq n}
\e (L_{i})^{d_{i}})[\MM_{g}^{n}]$$ is invariant under
permutations of the marked points, so $$m_{I}(T_{\Phi}^{*})
[\MM_{g}^{n}] = \frac{n!}{I!} \langle \tau^{I} \rangle \;,$$
where $I! = \prod_{k \geq 0} i_{k}!$. This rational number is the
image of $$\Phi_{g*}^{n}m_{I}(T_{\Phi}^{*}) \in
H^{*}(\MM_{g},\Q)$$ under the Gysin homomorphism defined by the
map from $\MM_{g}$ to a point; the latter is just evaluation
on the fundamental class of $\MM_{g}$, and the free energy [31 \S
2.16] can thus be recognized as the sum over $g$ and $n$ of terms of
the form $$\cl ([\Phi_{g}^{n}]/n!)[\MM_{g}] = \sum \langle 
\tau^{I} \rangle \frac{t^{I}}{I!} \;.$$
The classical Jacobian is similarly a configuration 
space of divisors, and the 
sum $$j_{g} = \sum_{n \geq 0} [\SP^{n}C_{g}]v^{n} \in MU^{0}_{\Q[v]}(\M_{g})$$
of symmetric powers of the modular family of curves is the pullback of
the product of the universal torus bundle in $MU^{-2g}(BSp(2g,\Z))$ by the 
coupling constant $$v^{g}\sum_{n \geq 0} [\C P(n)]v^{n} \; .$$ It is somewhat
surprising that no expression seems to be known for the class of 
this universal bundle.
 
\subsection*{1.4} The construction $\tau$ can be extended to 
stable curves which are not necessarily connected. From
the point of view of homotopy theory the space of
unordered configurations of points in $X$ is $$Q(X) 
= \coprod_{n \geq 0} E\Sigma_{n} \times_{\Sigma_{n}} X^{n}$$
but over the rationals this has the homotopy type of the
infinite symmetric product $\SP^{\infty}(X)$. We will therefore
take $$Q(\MM) = Q(\coprod_{g \geq 0} \MM_{g})$$ as a model for the 
space of stable curves, connected or not; its rational homology 
is the symmetric tensor algebra on $\oplus_{g \geq 0} H_{*}(\MM_{g},\Q)$, 
with $\MM_{0}$ and $\MM_{1}$ interpreted as one-point spaces. It is natural 
to think of the fundamental homology class of the orbifold $\SP^{n} \MM_{g}$ 
as the quotient in this tensor algebra of $[\MM_{g}]^{n}$ by $n!$, so 
$$[Q(\MM)] = \exp(\sum_{g \geq 0} [\MM_{g}]v^{3(g-1)})$$ defines a 
kind of fundamental class for $Q(\MM)$ with coefficients in $\Q[v]$.
Because no marked curve is stable when $3(g-1) + n < 0$, surfaces of 
small genus play a slightly anomalous role in these formulas; for 
example, it is useful to define $[\MM_{1}] = - \frac {1}{12}$.

The hard Lefschetz theorem (or the theory of mixed Hodge structures)
defines an action of $sl_{2}$ on the rational homology of a projective 
orbifold, with its fundamental class as 
a highest weight vector, so the homology of $Q(\MM)$ inherits an action
of $sl_{2}[v,v^{-1}]$ in which the class $[Q(\MM)]$ has 
conformal weight zero [18 \S 2.6].
%i.\e.\  $$(v \partial_{v} - H)[Q(\MM)] 
%= 0$$ with $H$ the usual diagonalizable basis vector in $sl_{2}$. 
The construction of $\tau_{g}$ extends multiplicatively to define an 
element $$\tau \in MU^{0}_{\Q[v]}(Q(\MM))$$ which sends this 
fundamental class to the partition function $$\tau_{W} 
= \tau_{*}[Q(\MM)] = \exp(\sum_{g \geq 0} F_{g}) \in 
\Q[t_{k}|k \geq 1][[t_{0},v]] \; .$$
%With the conventions of this paper, $t_{k}$ has homological
%degree $2k$, while the corresponding variables in Witten's formula
%are ungraded; but the terms in the expression above can be regrouped
%so that only terms of the form $v^{k}t_{k}$ appear. 

These results from physics suggest that the (immensely complicated) 
moduli space of curves has quite interesting homotopy-theoretic 
approximations, but (unlike the somewhat similar situation in algebraic 
$K$-theory) we do not yet understand these stabilizations in terms of 
universal mapping properties. That $\tau$ takes the fundamental class
of the moduli space to a highest-weight vector for the natural endomorphisms
of $\Delta_{\Q}$ is the central result of Kontsevich-Witten theory, but it is 
not yet a characterization.

\subsection*{1.5} The parameters $t_{k}$ have an intrinsic interpretation as
polynomial generators for the extended Landweber-Novikov Hopf algebra $$S = 
\Z[t_{k}|\; k \geq 0][t_{0}^{-1}]$$ of cooperations in complex cobordism. The universal 
stable cohomology operation is a ring homomorphism $$MU^{*}(X) \rightarrow 
MU^{*} \otimes S$$ and the characteristic number homomorphism of \S 1.1 can be 
defined as the composition of this map with the Thom map from cobordism to
ordinary cohomology. The universal operation sends the Euler 
class ${\e}$ of a complex line bundle to $$t({\e}) = \sum_{k \geq 0} 
t_{k}{\e}^{k+1} \; ;$$ it follows that the algebra $S$ represents the group 
of formal origin-preserving diffeomorphisms of the line, and that
the characteristic number homomorphism induces the classifying
map for the universal formal group law $$X +_{S} Y = t(t^{-1}(X) + 
t^{-1}(Y))$$ of additive type. The Lie algebra of formal vector fields
on the line is closely related to the Virasoro algebra, but the cobordism 
ring is more closely related to the untwisted charge one basic representation 
[21] than to the representation [31 \S 2.59] defined by topological gravity.

\section*{\bf \S 2 Background from algebra}

\subsection*{2.1} The simplest definition of the algebra $\Delta$ 
of $Q$-functions is as the quotient of the polynomial algebra
$\Z[q_{k}|k \geq 1]$ by the ideal generated by the coefficients of the relation
$$q(T^{\half})q(-T^{\half}) = 1 \; ,$$ where $$q(T^{\half}) = \sum_{k \geq 0} 
q_{k}T^{\half k}$$ is a generating function with $q_{0} = 1$ [14 \S 7]. This
algebra has a natural grading, which does not fit very comfortably with the
conventions of algebraic topology; we will assign the formal variable 
$T^{\half}$ cohomological degree one, but we will not assume that elements
of odd degree anticommute. The diagonal 
homomorphism $$q_{i} \mapsto \sum_{i=j+k} q_{j} \otimes q_{k}$$ makes $\Delta$
into a bicommutative Hopf algebra over $\Z$, and the relation $$(-1)^{i}q_{i}^{2} 
= 2q_{2i} + 2\sum_{i-1 \geq k \geq 1} (-1)^{k-1} q_{k}q_{2i-k}$$
implies that the square of a generator can be expressed as a sum of 
monomials in which no $q_{k}$ appears with exponent greater than one. It
follows that $\Delta$ is a free module over the integers, with a 
basis of square-free monomials; similarly, $\Delta[\half]$
is a polynomial algebra generated by the elements $q_{2k+1}$, while the reduction
of $\Delta$ modulo two is an exterior algebra. Being torsion-free, 
$\Delta$ embeds in $\Delta_{\Q} = \Delta \otimes \Q$,
and its defining relations imply that the power series $$\log q(T^{\half}) 
= \sum_{i \geq 0} \frac{2x_{k}}{2k+1} T^{k + \half}$$ 
contains only odd powers of $T^{\half}$. Newton's identity $$(2k+1)q_{2k+1} = 
2(x_{0}q_{2k} + x_{1}q_{2k-2} + \dots + x_{k})$$ shows that the classes 
$2x_{k}$ are integral, e.g.\ $2x_{0} = q_{1}$ and $2x_{1} = 3q_{3} - 
q_{1}q_{2}$.  

\subsection*{2.2} The generators $q_{i}$ of $\Delta$ can be interpreted as specializations at 
$t = - 1$ of Hall-Littlewood symmetric functions $q_{k}(\Lambda;t)$ of 
the eigenvalues of a positive-definite self-adjoint matrix $\Lambda$, defined by
$$q_{\Lambda,t}(T^{\half}) = \sum q_{k}(\Lambda;t) T^{\half k}  
=  \det \; \frac{1 - \Lambda^{-1}T^{\half}t}{1 - \Lambda^{-1} T^{\half}} 
\; ;$$ the primitive element $x_{k}$ thus becomes 
the power sum $tr\; \Lambda^{-2k-1}$. The canonical symmetric bilinear 
$\Z[t]$-valued form on the algebra of Hall-Littlewood functions defines a 
positive-definite inner product on $\Delta_{\Q}$ regarded as
a quotient of that algebra [23 III \S 8.12]. The product and coproduct of $\Delta$ are 
dual with respect to this bilinear form; in fact $\Delta [\half]$ is a 
positive self-adjoint Hopf algebra in the sense of Zelevinsky 
[23 I \S 5 ex 25]. The classical $Q$-functions are the orthogonalization 
of the basis of square-free monomials with respect to this inner product.   

\subsection*{2.3} The inner product on $\Delta_{\Q}$ defines a skew bilinear
form on the complexification of its
space of primitives, which allows us to interpret $\Delta_{\Q}$
as a standard representation of a Heisenberg algebra, or
alternately [26, appendix I] of the group of antiperiodic loops on the circle.
The latter group has two connected components, and in
some ways [13, appendix 8] it is natural to think of $\Delta$ as
$\Z \times \Z/2\Z$-graded algebra, with $q_{k}$ in homological
degree $(k-1,1)$; but in other contexts the $\Z/2\Z$-grading appears
as an action of the Galois group of $\Q(\sqrt 2)$ over $\Q$. 

If, using the notation of [32 \S 1.7], we write $$2^{\half}
\alpha_{-k-\half} = - x_{k}$$ when $k > 0$ then the operators
$$L_{0} = \sum_{k>0} \alpha_{-k} \alpha_{k} + \frac{1}{16}\;,\;
L_{n} = \half \sum \alpha_{k} \alpha_{n-k} \;,\;n \neq 0$$ define
a twisted [9 \S 9.4] charge one representation of the Virasoro
algebra on $\Delta_{\Q}$. There appear to be deep general connections
between self-dual Hopf algebras and vertex operator algebras [1,13]. 

\subsection*{2.4} It is striking that these $Q$-functions have been important in 
algebraic topology for more than a generation: in Cartan's 1960 
Seminar [4 \S 17] the integral homology of $Sp(\HH)/\U = \Omega^{5} 
\mathbb O$ appears as the quotient $\Delta^{+}$ of a polynomial algebra 
on generators $q^{+}_{k}$ of degree $2k$, modulo the ideal generated
by the coefficients of the relation $$q^{+}(T)q^{+}(-T) = 1 \;
,$$ where $$q^{+}(T) = \sum_{k \geq 0} q_{k}T^{k} \; ;$$
the algebras $\Delta$ and $\Delta^{+}$ thus differ only by a
doubling of the grading. The cohomology of the corresponding finite-dimensional
homogeneous spaces has been studied more recently along similar lines [16].

In light of the results of Karoubi cited in \S 1.2, an algebra
of $Q$-functions forms a substantial part of the stable integral homology 
of the space of Abelian varieties; the primitives $x_{k}^{+}$ are dual
to Mumford's classes $\lambda_{2k+1}$. The stable rational cohomology of the 
Riemann moduli space contains a polynomial algebra on the Mumford classes 
$\kappa_{k}$ of degree $2k$, and the dual homology Hopf algebra
contains the polynomial algebra [24] on dual primitives $\hat
\kappa_{k}$. The map of \S 1.2 kills the even classes $\hat
\kappa_{2k}$, while $$\lambda_{2k-1} \mapsto (-1)^{k} \frac {B_{k}}{2k}
\kappa_{2k-1} \; ,$$ where $B_{k}$ is the $k$th
Bernoulli number [29 \S 6.2]; the even Mumford classes are
detected by the constructions of [28]. Stabilization implies
that $$\hat \kappa_{k} \in H_{2k}(\M_{g},\Q)$$ for
sufficiently large $g$, and (since these classes vanish on
decomposables) it follows from the characteristic number formula
of \S 1.1 that $$\tau_{g*}(\hat \kappa_{k}) = vt_{k+1} \; ;$$
the maps $\tau_{g}$ are however not ring homomorphisms with respect
to the usual multiplication [26] on complex cobordism.

\subsection*{2.5} The importance of the theory of symmetric
functions in the work of Kontsevich and Witten was discovered by
Di Francesco, Itzykson, and Zuber [7 \S 3.2], but it was J\'ozefiak [17]
who saw the connection with $Q$-functions. In our formalism, 
their map $$t_{k} \mapsto
- (2k-1)!! \; tr \; \Lambda^{-2k-1} = - (2k-1)!! x_{k} \; ,$$
[where the `odd' factorial $(2k+1)!!$ is the product of the odd 
integers less than or equal to $2k+1$, with $(-1)!! = 1$ by convention] 
defines a homomorphism from the extended
Landweber-Novikov algebra to $\Delta[\half]$. The image of the partition
function $\tau_{W}$ under this map satisfies a large family of differential
equations [22] which can be summarized impressionistically by the assertion that
$$\tilde \tau_{W}(x_{i}) = \tau_{W}(x_{i} + \delta_{i,1})$$ is an 
$sl_{2}$-invariant highest weight vector for the natural Virasoro action 
on (a completion of) $\Delta_{\Q}$. The formal series $\tau_{W}$ is 
divergent at $x_{i} = \delta_{i,1}$, but this claim can be reformulated 
precisely as $$\tilde L_{n} \tau_{W} = 0 \; , \; n \geq -1 \; ,$$ where 
$\tilde L_{n}$ is the linear shift of $L_{n}$ defined by the map which 
sends $x_{1}$ to $x_{1} - 1$, cf.\ [19]. This defines a charge one vertex
operator algebra embedded in a completion of $\Delta_{\Q}$, cf.\ [8].

There is reason [22 \S 7, 23 III \S 7 ex 7] to expect that these results 
generalize to Hall-Littlewood functions at other roots of unity.

\section*{\bf \S 3 The Kontsevich-Witten genus}

\subsection*{3.1} Quillen's theorem establishes a bijection 
between one-dimensional formal group laws over a commutative ring $A$, and 
homomorphisms from the complex cobordism ring to $A$, so we can define
a formal group law $$X +_{Q} Y = \kw^{-1}(\kw(X) + \kw(Y))$$ over the 
localization $\Delta[q_{1}^{-1}]$, and hence a $Q$-function valued genus of
complex manifolds, by specifying its exponential series to be
$$\kw^{-1}(T) =  \sum_{k \geq 0} (2k-1)!! x_{0}^{-1}x_{k} 
T^{k+1} \; .$$ The homomorphism $$\kw : MU \rightarrow
S \rightarrow \Delta[q_{1}^{-1}]$$ classifying this group law factors through 
the map classifying the universal additive law of \S 1.5 and therefore 
defines a group law of additive type. 

\subsection*{3.2} This Kontsevich-Witten genus is defined by a 
nonstandard orientation on an ordinary cohomology theory represented 
by the generalized Eilenberg-Mac Lane space associated to 
$\Delta[q_{1}^{-1}]$; its vertex algebra structure is the source of 
Witten's Virasoro operations. 

Some properties of this genus are unfamiliar to the point of 
pathology: not only is its exponential series integral, for example, but 
modulo an odd prime it is polynomial as well. The point of this section is that 
the orientation defining this genus can be interpreted as the formal completion of a coordinate, 
in the sense of [10 \S 1.8], on a $\T$-equivariant cohomology 
theory taking values in sheaves of modules over a certain abelian group 
object in the category of ringed spaces. This global object is 
semi-classical, and some of the strangeness of the Kontsevich-Witten genus
is a property not of the group object itself, but of a rather inconvenient
coordinate on it.  

It seems simplest to present this 
construction in two parts: the first step will define a local version 
over a basic group object $\X$. The genus itself will then be defined 
by a family of group objects parametrized by a Grassmannian of 
positive-definite self-adjoint matrices.

\subsection*{3.3.1} The basic idea comes from the classical theory 
of functions of one complex variable: if $g$ is an entire function,
of exponential type in every right half-plane, then under
certain circumstances there is an asymptotic relation
$$\sum_{n \in \Z} g(n)z^{n} \sim 0 \; ;$$
such equations may look more familiar written in the form
$$ \sum_{n \geq 0} g(n)z^{n} \sim - \sum_{n \geq 1} g(-n)z^{-n}
\; .$$ If $G(z)$ is the left-hand sum, $\check G(z^{-1})$ will denote 
the sum on the right.
When $g$ is rational [e.g.\ $g = 1$] such relations
are familiar, but a less trivial example is defined by
$$g(w) = \Gamma(1 + \alpha w)^{-1}$$ with $0 < \alpha < 2$. The
resulting (Mittag-Leffler) function $$\exp_{\alpha}(z) = \sum_{n
\geq 0} \frac{z^{n}}{\Gamma(1 + \alpha n)}$$ thus has the asymptotic
expansion $$ \exp_{\alpha}(z) \sim - \sum_{n \geq 1} \frac{z^{-n}}{\Gamma(1
- \alpha n)}$$ for $z$ outside a sector of angular width $\alpha
\pi$ centered on the positive real axis, cf. [12 \S 11.3.23]. 
This function is especially interesting when $\alpha$ is a
rational number between zero and one, but we will
be concerned only with $\exp_{\half}(z)$ which, up to 
normalization, is the Laplace transform of Gaussian measure. 
%It is also a solution of the differential equation 
%$$y'' - xy' - y = 0 \; ,$$ and moreover
%$$\exp_{\half}(x) = e^{x^{2}} [ 1 + E(x) ]$$
%for $x$ real, where $$E(x) = 2 \pi^{- \half} \int^{x}_{0}
%e^{-t^{2}} dt$$ is the normalized Gaussian error integral. 
Using the duplication formula for the gamma function, its asymptotic
expansion takes the form $$\exp_{\half}(z) \sim - \pi^{-
\half}z^{-1} \sum_{n \geq 0}(2n-1)!! (-2z^{2})^{-n}$$ for $z$
outside a sector of width $\half \pi$ centered on the positive
real axis; in particular the expansion is valid along the entire
imaginary axis. The odd function $$\sin_{\half}(z) =
- \frac{i}{2} [\exp_{\half}(iz) - \exp_{\half}(-iz)]$$ therefore
satisfies $$\sin_{\half}(x) \sim \pi^{-\half}x^{-1} \sum_{n \geq
0} (2n-1)!! (2x^{2})^{-n}$$ as $x$ approaches infinity in
either direction along the real axis.

\subsection*{3.3.2} Now consider the behavior on the real line 
of the entire function
$$\epsilon(z) = (\frac {\pi}{2z})^{\half} \sin_{\half}
((\frac {z}{2})^{\half})$$ defined by the power series $$\frac
{\pi^{\half}}{2} \sum_{n \geq 0} \frac {(-\half z)^{n}}{\Gamma(n
+ \frac{3}{2})} = \sum_{n \geq 0} \frac{(-z)^{n}}{(2n+1)!!} \; .$$ 
The argument above implies that 
$$\epsilon(x) \sim \check \epsilon(x^{-1})
= \sum_{n \geq 0} (2n-1)!! x^{-n-1}$$ for $x$ large and positive;
it follows that $\epsilon$ is monotone decreasing for
positive $x$. On the other hand it is clear from its power
series expansion that $\epsilon$ is monotone decreasing on the
negative real axis, so $$\epsilon : (\R,0) \rightarrow
(\R_{+},1)$$ is a bijection. The open 
interval $$\X = (0,\infty)^{+} - \{ 1 \}$$
can now be made an abelian group by interpreting $\epsilon$ to be the exponential
of a group law defined in a neighborhood of infinity on the projective 
line: if $x$ and $y$ are large (i.\ e.\ nonzero) real numbers, then 
$$x +_{\infty} y = \frac {xy}{x+y}$$ will again be a large real
number. It is easy to see that the resulting composition 
on $P^{1}(\R) - \{0\}$ is
associative, with $\infty$ as identity element; the inverse of
$x$ is just $-x$. We can therefore construct a formal group
law $+_{\infty \epsilon}$ over $\R$ by requiring that $$\epsilon(x)
+_{\infty \epsilon} \epsilon(y) \sim \epsilon(x +_{\infty} y)$$
for $x$ and $y$ large and positive. This group law as the completion 
at the identity of the analytic composition $$X +_{\epsilon} Y = 
\epsilon (\epsilon^{-1}(X) +_{\infty} \epsilon^{-1}(Y))$$ with 
translation-invariant one-form $d\epsilon^{-1}(T)^{-1}$; 
its exponential is the series $\check \epsilon$. The involution 
$$[-1]_{\epsilon}(T) = \epsilon (- \epsilon^{-1}(T))$$ interchanges 
$(0,1)$ and $(1,\infty)$, identifying $\X$ with the group
completion of the semigroup defined on $(0,1)$ by $+_{\epsilon}$.
This group object is smooth, but not analytic; its exponential
map has trivial domain of convergence.

\subsection*{3.3.3} The usual sheaf of smooth real-valued functions defines 
the structure of a ringed space on $\X$, but its law of addition is also 
compatible with the sheaf of real-analytic functions formally completed at 
the origin. The composition $$x \mapsto \Psi(x) = \exp( - \epsilon^{-1}
(x)^{-1}) : \X \rightarrow \R_{+}^{\times}$$ is a homomorphism to the real 
multiplicative group; pulling the sheaf over ${\mathbb G}_{m}(\R)$ defined 
by $\T$-equivariant $K$-theory back along this map yields a $\T$-equivariant 
cohomology theory taking values in the category of
sheaves of modules over $\X$. This is a kind of ordinary equivariant
cohomology theory with coefficients in $\X$; its Chern-Dold character sends the
standard one-dimensional representation of the circle to $\Psi$, regarded
as a section of the structure sheaf. The Hirzebruch genus $$\frac {x}{\epsilon
(x^{-1})} = \pi^{-\half} \frac {(2x)^{\half}}{\sin_{\half}((2x)^{-\half})}$$
of the associated complex orientation is the reciprocal of a power series in
$x^{-1}$ with trivial constant term. It defines a function on $P^{1}(\R)$
analytic aside from a jump discontinuity at 0.

\subsection*{3.4} By replacing the function $\epsilon$ with 
$$\epsilon_{\Lambda}(x) = (\frac {\pi}{2x})^{\half}
\; tr \; \sin_{\half} ((\frac {x}{2})^{\half} \Lambda )$$ we obtain
a $\T$-equivariant theory taking values in sheaves over a family of 
abelian group objects $\X_{\Lambda}$ parametrized by equivalence classes of
self-adjoint positive-definite matrices $\Lambda$. It is remarkable that 
if $\Lambda$ is an $np \times np$ matrix and $t$ is a $p$th root of unity, 
then it follows from the definition of \S 2.2 that 
$$q_{\Lambda,t}(T^{\half}) = q_{\Lambda^{-1},t^{-1}}(T^{-\half}) \; ; $$ in 
particular, if $\Lambda$ is an endomorphism of an even-dimensional vector 
space, and $\check q(T^{\half})$ denotes the generating function 
for the Hall-Littlewood functions at $-1$ associated to the matrix 
$\Lambda^{-1}$, then $$\check q(T^{-\half}) = q(T^{\half}) \; .$$ It 
seems reasonable to expect that this $\T$-equivariant theory is a real version of
a theory over the quotient of some completion of $\Delta \otimes \check \Delta$ 
by the coefficients of the relation $$\check q(T^{-\half})q(-T^{\half}) = 1 \;.$$

\section*{Appendix: vertex operator algebras and self-dual Hopf
algebras}

This appendix has been added in July 1998 to the published 
version of this paper, which has appeared in `Homotopy theory 
via algebraic geometry and representation theory', ed. S. Priddy
and M. Mahowald, in Contemporary Mathematics. It is a kind of 
commentary on Andy Baker's construction of vertex operator algebras
associated to even unimodular lattices using ideas from the theory 
of (bicommutative) Hopf algebras. I suspect that these methods will
have further use, and I have tried here to summarize Baker's
results in the language of group-valued functors, and to include
some references to related work [e.\ g.\ on $\lambda$-rings [3]]
which might otherwise be overlooked.

\subsection*{A.1} The functor which assigns to a commutative ring $A$,
the abelian group (under multiplication) of formal series $$h(T)
= 1 + \sum_{i \geq 1} h_{i} T^{i} \in (1 + TA[[T]])^{\times} :=
\W_{0}(A) \; ,$$ is represented by the polynomial Hopf algebra $S
:= \Z[h_{i}| i \geq 1]$ with comultiplication $$\Delta (h_{i}) =
\sum_{i=j+k} h_{j} \otimes h_{k} \; .$$ $\W^{0}(A)$ is
essentially the classical Witt ring of $A$ [6 V \S 2, 2], and its
representing algebra can be identified with the usual ring $S$ of
symmetric functions. More precisely, there is a natural isomorphism 
of the set $\W_{0}(A)$ with the set of ring homomorphisms from $S$ to
$A$, such that the map induced by $\Delta$ agrees with
multiplication of power series. The functor $\W(A)$ which assigns
to $A$ the set of {\bf all} invertible power series over $A$ is 
represented by the tensor product $S[h_{0},h_{0}^{-1}]$ of the 
usual ring of symmetric functions with a Hopf algebra
representing the multiplicative groupscheme. It will be useful 
to know that $\W_{0}(A)$ has a
natural commutative ring-structure $*$ characterized by the
identity $$(1 + aT)*(1+bT) = (1 + abT) \; .$$

If $\mL$ is a free abelian group of finite rank $l$
[e.\ g.\ a lattice, with dual $\mL^{*} = \Hom_{\rm ab}(\mL,\Z)$],
then it is easy to see [for example by choosing a basis] that the
functor $$A \mapsto \Hom_{\rm ab}(\mL,\W_{0}(A)) = \mL^{*}
\otimes_{\Z} \W_{0}(A)$$ is also represented by a Hopf algebra,
which can be taken to be an $l$-fold tensor product of copies of
$S$. I will write $\otimes^{\mL^{*}} S$ for this representing
object, which Baker calls $S(\mL)$. The point of the definition
at the beginning of [1 \S 3] is to present a coordinate-free
definition of this object: a homomorphism from $\mL$ to
$\W_{0}(A)$ defines a family $\lambda \mapsto h_{i}(\lambda)$ of
maps from $\mL$ to $A$, such that $$h_{i}(\lambda + \lambda') =
\sum_{i=j+k}h_{j}(\lambda) h_{k}(\lambda') \; ;$$ the ring
$\otimes^{\mL^{*}}S$ thus represents the tensor product functor
in a natural way, without specifying a basis. The generating
function $$h^{\lambda}(T) = \sum_{i \geq 0} h_{i}(\lambda)
T^{i}$$ is a convenient substitute for such a choice.

It is tempting to think of $\mL^{*}$ as a constant groupscheme,
and to interpret $\mL^{*} \otimes_{\Z} \W_{0}$ as a tensor
product in a category of group-valued functors. A natural
internal product on a suitable category of commutative and
cocommutative Hopf algebras has been studied by Goerss
[11] and by Hunton and Turner [14]. 

In fact $\otimes^{\mL^{*}}S$ is only a part of the vertex
operator algebra associated to a lattice; the full construction
is usually described as the graded tensor product $\Z[\mL]
\otimes (\otimes^{\mL^{*}}S)$ of the symmetric functions with the
group ring of $\mL$. In this context it is natural to think of
this group ring as graded, with the element $\lambda$ in degree
$\la \lambda,\lambda \ra$, and to give $S$ its usual grading. In
view of the discussion above we can define $\V(\mL)$ to be the
Hopf algebra representing the functor $\mL^{*} \otimes_{\Z} \W$.

\subsection*{A.2} If $H$ is a commutative and cocommutative Hopf algebra,
projective of finite rank over a base ring $k$, then the module
$$H^{*} = \Hom_{k}(H,k)$$ is again a bicommutative
Hopf algebra. If $$A \mapsto \Hom_{k-{\rm alg}}(H,A) := \HH(A)$$
is the group-valued functor $H$ represents, then its dual Hopf
algebra $H^{*}$ represents the Cartier dual functor $$A \mapsto
\Hom_{\rm gp-valued\;functors}(\HH,\G_{m})(A) \; ,$$ cf.\ [6 II \S
1 no.\ 2.10]. A self-duality on such a (commutative and
cocommutative, projective and finite) Hopf algebra can thus be
defined either as an isomorphism $$H \rightarrow H^{*}$$ of Hopf
algebras, or as an isomorphism $$\HH \rightarrow \HH^{*}$$ of
abelian group-valued functors; alternately, such a structure can
be defined either as a nondegenerate pairing $$H \otimes H
\rightarrow k$$ of algebras [with suitable properties [23 I \S 5
ex.\ 25, 13]] or as a nondegenerate pairing $$\HH \times \HH
\rightarrow \G_{m}$$ of abelian group-valued functors.

Similar dualities exist more generally, in particular in the
context of locally finite graded Hopf algebras. I will write
$S_{*}$ for the ring of symmetric functions given its usual
grading, and $S^{*}$ for its graded dual; the Hall inner product
discussed by MacDonald then defines an isomorphism $$S^{*}
\rightarrow \check S_{*}$$ of graded algebras, with the
convention that $\check H_{*}$ is the graded algebra $H_{-*}$,
i.\ e.\ with its grading negated. A self-duality on a locally
finite graded Hopf algebra $H$ can thus be defined by a
nondegenerate graded pairing $$H \otimes \check H \rightarrow k$$
with suitable properties. If we think of the grading on $H$ as an
action of the multiplicative group on $\HH$, then $\check \HH$
will have the inverse action.

In terms of group-valued functors, the Hall duality
defines a homomorphism $$\W_{0} \times \W_{0} \rightarrow
\G_{m}$$ which is familiar classically as part of the theory of
the Artin-Hasse exponential [6 V \S 4 no. 4.3]. It seems to be
simpler to work not with graded rings but with the formal
completion of $\W_{0}$ at the origin, i.\ e.\ to consider the
subfunctor $\hat \W_{0}$ defined on complete local rings $A$ by
series $h(T)$ with coefficients $h_{i}$ in the maximal 
ideal of $A$ such that $h_{i} \rightarrow 0$ as $i \rightarrow
\infty$. Since $\W_{0}$ is a commutative ring-valued functor, the
duality can be constructed in terms familiar from the theory of
Frobenius algebras. It will suffice to define a suitable `trace'
morphism from $\hat \W_{0}$ to $\G_{m}$; the pairing will then be
the result of this trace applied to the Witt ring product. But 
if $h \in \hat \W_{0}(A)$ then $h(1) \in A^{\times}$, and $$h,g 
\mapsto (h*g)(1)$$ will do the trick \dots

\subsection*{A.3} The principal observation in this appendix 
is that some of
Baker's formulae can be simplified considerably by use of the
self-dual nature of the underlying Hopf algebra. His formula 2.4
can be summarized as follows: let $$\Y(h(w)) = h(z + w)
\otimes \check h(-(z + w)^{-1}) \in (S \otimes \check
S)[[z,z^{-1}]][[w]]$$ be defined by expanding the coefficient
in the right-hand term as $\sum_{i \geq 0} (-w)^{i}z^{-i-1}$; 
here $\check h(T) \in \check S[[T]]$ is the analogue of the 
generating function $h(T)$. This formula extends
multiplicatively to define a homomorphism $$\Y : S \rightarrow (S
\otimes \check S)[[z,z^{-1}]]$$ of (graded) ring{\bf oids}, the
point being that multiplication is not always defined in the
object on the right. The assertion, more precisely, is that for
any $b$ and $b'$ in $S$, the product $\Y(b)\Y(b')$ is in fact
defined, and equal to $\Y(bb')$. This is not quite Baker's
definition of the vertex operator, but the two agree under the
identification of $S \otimes \check S$ with $\End(S)$ given by
Hall duality: the commutative product in the former algebra is a
resource not available in the latter, and it is just what
is needed to make sense of the assertion that Baker's formula is
multiplicative. 

The definition of Cartier duality in terms of an
internal $\Hom$ functor on the category of commutative
groupschemes makes it natural to identify $\mL \otimes_{\Z}
\W_{0}$ with $\Hom(\mL^{*} \otimes_{\Z} \W_{0} ,\G_{m})$; the
isomorphism $$\mL^{*} \rightarrow \mL$$ defined by the pairing on
a self-dual lattice thus extends to a (grade-negating)
isomorphism between $\otimes^{\mL^{*}}S$ and $\otimes^{\mL}
\check S$. The analogue of the generating function
$h^{\lambda}(T)$ is a generating function $\check h^{\lambda}(T)
\in (\otimes^{\mL} \check S)[[T]]$, and the formula above for
$\Y$ extends immediately to this more general context [1 \S 3.3].
To define $\Y$ on the whole of $\V(\mL)$, it remains to construct
$\Y$ on elements of the group ring of $\mL$ \;: Baker's formula is
$$\Y(\lambda) = \lambda h^{\lambda}(z) \otimes \check
h^{\lambda}(-z^{-1}) z^{\lambda} \; ,$$ where $z^{\lambda} \in
\End(\Z[\mL])[z,z^{-1}]$ is the operation $\mu \mapsto \mu z^{\la
\lambda,\mu \ra}$ \dots \bigskip

Much of the preceding construction appears to
generalize in a relatively straightforward way from the classical
algebra of symmetric functions to the Hopf $\Z[t]$-algebra $HL$
of Hall-Littlewood symmetric functions [23 III \S 5 ex.\ 8]: this
can be defined as the polynomial algebra on generators $q_{i}$
related to the power sum symmetric functions $p_{n}$ by the
formula $$q(T) = \sum_{i \geq 0} q_{i} T^{i} = \exp (\sum_{n \geq
1} ( 1 - t^{n}) p_{n} \frac{T^{n}}{n})$$ [23 III \S 2.10]. As $t
\rightarrow 0, \; q(T) \rightarrow h(T)$, which suggests that Baker's
formula for $\Y$ might have an interesting analogue in $HL$.
Taking formal logarithms, we can rewrite that formula in terms of
`reduced' power sums $\tilde p_{n} = p_{n}/n$ as $$\Y(\tilde p_{n}) 
= \sum_{m \in \Z - \{0\}} \frac{1 - t^{|m|}}{1 - t^{n}} \binom{m}{n} 
\tilde p_{m} z^{m-n} \;,$$ where $n$ is positive and $p_{-m} = (-1)^{m+1} 
\check p_{m}$ for $m$ negative. 

The theory of $Q$-functions is one of the original
motivations for this account: these are elements of the self-dual
Hopf algebra which represents the functor $$A \mapsto \{ h \in
\W_{0}(A) \; | \; h(-T) = h(T)^{-1} \} \; .$$ An analogous
construction exists for any prime $p$: if $\omega$ is a primitive
$p$th root of unity, the subfunctor $$\{ h \in \W_{0} \; | \;
\prod_{0 \leq k \leq p-1} h(\omega^{k}T) = 1 \}$$ is also
represented by a self-dual Hopf algebra. Expressed in terms of
power sums, the relation defining this quotient becomes $p_{n} =
0$ when $p$ divides $n$. It is also possible to think of these
algebras as quotients of $HL$ defined by specializing $t$ to a
root of unity [23 III \S 7 ex.\ 7]; perhaps requiring that its order
be prime is superfluous. Considerable effort [13,15] has been
devoted to constructing aspects of a vertex algebra structure on
$HL$ and on the rings of $Q$-functions, and it may be worth
noting that when $t$ is a primitive $p$th root of unity, the
formula above for $\Y$ continues to make sense on the quotient of
$HL$ defined by sending $p_{n}$ to zero when $p$ divides $n$; but I
am not sure if this fits with the thinking in [9 \S 9.2.51], and 
at the moment I don't know how to fit a lattice into the picture. 

As Baker notes, some constructions of VOA theory are
strikingly close to constructions familiar in algebraic topology.
As a closing footnote I'd like to suggest the possible relevance
of the classifying space $BU_{\T}$ for $\T$-equivariant
$K$-theory in this context: the vertex operator $\Y$ can be
interpreted as a homotopy-class of maps $$BU_{\T} \times BU_{\T}
\rightarrow BU_{\T}[[z,z^{-1}]] \; ,$$ or in terms of
representation theory as a construction which relates [or
`fuses'] two representations of $\T$ specified at the points $0,1$ 
on the projective line, to define a representation at $\infty$. This 
seems very close to the point of view of Huang-Lepowsky [and Segal] \dots

\bibliographystyle{amsplain}

\end{document}